\documentclass[a4paper,leqno,12pt]{article}
\usepackage{amsmath,amsfonts,amssymb}
\usepackage{fancyhdr}
\usepackage{a4wide}
\usepackage{hyperref}
\usepackage{makeidx}
\makeindex 
 
\newcounter{EQNR}[NNN] 
\newcommand{\nne}{\refstepcounter{EQNR}  \tag{\theNNN .\theEQNR} }
\newcommand{\nn}{\vskip 2mm 
\refstepcounter{NNN}\noindent {\bf \theNNN . }}

\nopagebreak

 \DeclareMathOperator{\Ei}{Ei}

\newcommand{\ta}{{}^t\!}

\newcommand{\SL}{{\rm SL}}

\newcommand{\R}{{\bf R}}

\newcommand{\Z}{{\bf Z}}
\newcommand{\N}{{\bf N}}

\newcommand{\Ha}{{\mathfrak H}}

\newcommand{\SO}{{\rm SO}}

\newcommand{\Ima}{\,{\rm Im}\,\,}
\newcommand{\Rea}{\, {\rm Re}\,\,}

\begin{document}

\thispagestyle{empty}

{\bf \Large On Kudla's Green function for signature (2,2) Part II}\\[1cm] 
{\bf By Rolf Berndt and Ulf K\"uhn}\\[.5cm]

\begin{abstract} Around 2000 Kudla presented conjectures about deep relations between
arithmetic intersection theory, Eisenstein series and their derivatives, and special values of Rankin $L-$series.
The aim of this text is to work out the details of an
old unpublished draft on the second  author's attempt to prove these conjectures for the
case of the product of two modular curves. 

 In part one we proved
that the generating series of certain modified arithmetic special cycles is as predicted by Kudla's conjectures 
a modular form with values in the first arithmetic Chow group. 
Here we pair this generating series with the square of the first arithmetic Chern class of the line bundle of modular forms. Up to previously known Faltings heights 
of Hecke correspondences only 
integrals of the Green functions $\Xi (m)$ over $X$ had to be computed.
The resulting arithmetic intersection numbers turn out to be as predicted by Kudla to be strongly related to the Fourier coefficients of the derivative of the classical real analytic Eisenstein series $E_2(\tau ,s)$.   

\end{abstract}

\tableofcontents

\section{Introduction}

As in the first part of our work \cite{BK1} we
  consider the natural models of 
  product of modular curves $X=X(1)\times X(1)$ 
and its Hecke correspondence $T(N)$ over the integers $\mathbb Z$.
We had introduced the modified arithmetic special cycles
\begin{align*}
\widehat{Z}_\rho (m) := (T(m),\widetilde{\Xi }_\rho (v,z,m)) \in \widehat{CH}^1(X)  
\end{align*}
 and a modified Kudla generating series
\begin{align}
\widehat \phi_{K,\rho}= \sum \widehat{Z}_\rho (m)q^m .\nne
\end{align}
We proved in \cite{BK1}that $\widehat \phi_{K,\rho}$ is a modular form for $\SL_2(\Z)$ of weight $k$ with coefficients in the arithmetic Chow group $\widehat{CH}^1(X)$.
 
Hence, for all linear maps $L: \widehat{CH}^1(X)  \rightarrow \mathbb R$ the series 
$L(\widehat{Z}_\rho (m))q^m$ is a (nonholomorphic) $\mathbb R-$valued modular form in the usual sense. 
Now we denote by $\widehat{c}_1( \overline{\mathcal L})$ the first arithmetic Chern class\footnote{Note this is an arithmetic cycle in the arithmetic 
Chow group with loglog-growth in the sense of \cite{BKK}}
of the line bundle of modular forms $\overline{\mathcal L}(12,12)$ of bi-weight $(12,12)$ equipped with the Petersson metric. Then we
choose the linear map $L(-) = \widehat{c}_1( \overline{\mathcal L})^2\cdot (-)$ and prove 
the following result, which complements the work of Kudla, Rapoport and Yang 
in the $O(1,2)$ case \cite{KRY}, \cite{KRY1} and provides
a first confirmation of Kudla's conjectures \cite{Ku} in dimension $2$.\\

{\bf Main Theorem} (modified Kudla conjecture){\bf.}
\emph{We have an identity of modular forms
  \begin{align}\label{eq:kudla-modi}
\widehat{c}_1( \overline{\mathcal L})^2\cdot \widehat \phi_{K,\rho}
 = \,\mathbb E'_2(\tau ,1) + f_\rho (\tau )\nne
\end{align}
with $\mathbb E'_2(\tau ,1)$  the derivative  of a non-holomorphic Eisenstein series $\mathbb E_2(\tau ,s)$ with respect to $s\in\mathbb C$  and a certain modular form $f_\rho (\tau )$.}\\

The above Eisenstein series equals
\begin{align*}
\mathbb E_2(\tau ,s) := -12\psi (s)E_2(\tau ,s),  
\end{align*}
where $\psi $ is a meromorphic function with
\begin{align*}
\psi (s) = -1 + 4(\frac{\zeta '(-1)}{\zeta (-1)} + \frac{1}{2})(s-1) + O((s-1)^2)  
\end{align*}
and the expansion of $E_2(\tau,s) = (1/(2\pi i))\partial _\tau E^\ast (\tau ,s)$ (c.f. \eqref{Estar}) is  determined  as follows.\\

{\bf Theorem.} \label{ThEiss}  \emph{The coefficients of  the Fourier expansion  of the weight $2$ Eisenstein series for $\SL_2(\Z)$
\begin{align*}
E_2(\tau ,1) = \sum_{m\in\Z}a(v,1,m) q^m\nne
\end{align*}
and those of the the derivative of $E_2(\tau ,s)$ with respect to $s$ 
\begin{align*} 
E'_2(\tau ,1) = \sum_{m\in\Z}a'(v,1,m) q^m
\end{align*}
are given by
\begin{align*} 
a(v,1,m) &= 
\begin{cases}
\sigma (m) = \sum_{d|m} d\,\,&{\rm for} \,\,m>0\\ 
0\,\,&{\rm for}\,\,m<0
\end{cases}\\
\intertext{and by}
a'(v,1,m) &= \begin{cases}
\sigma (m)(1/(4\pi mv) + \sigma '(m)/\sigma /m)) \,\,&{\rm for} \,\,m>0\\
\sigma (m)(\Ei(-4\pi |m|v) + 1/(4\pi |m|v)e^{-4\pi |m|v})\,\,&{\rm for}\,\,m<0
\end{cases}
\end{align*}
where with 
\begin{align*} 
\sigma ^\ast _{s}(m) := |m|^s\sum_{d|m}d^{-2s} 
 \end{align*}
we abbreviate $\sigma '(m)/\sigma (m):= {\sigma ^\ast }'_{1/2}(m)/\sigma ^\ast _{1/2}(m).$ 
}\\

 We also will calculate the constant terms 
$a(v,1,0)$ and $a'(v,1,0)$ below (c.f. Theorem \ref{ThderEiss})
although we  won't need them for this work

\nn {\bf Remark.}\label{RemmodEis} The Fourier expansion 
$E_2(\tau ,s) = \sum a(v,s,m)q^m $
translates via the multiplication by $\psi$ to 
$
\mathbb E_2(\tau ,s) = \sum A(v,s,m)q^m $,
where the first terms of the Taylor expansion at $s=1$ of the Fourier coefficients 
are given by
\begin{align*}\label{eqmodEis}
 A(v,1,m) &= 12 a(v,1,m),\\\nonumber
 A'(v,1,m) &=\, -48(\frac{\zeta '(-1)}{\zeta (-1)} + \frac{1}{2})a(v,1,m) +12 a'(v,1,m).\nne
\end{align*}\\

The main steps in the proof of our main theorem are to calculate
 and to compare both sides of \eqref{eq:kudla-modi} termwise. \\

Arakelov theory (\cite{BKK} Proposition 7.56) gives us for $m\not=0$   
  the relation
\begin{align}\label{cxcxZm}
\widehat{c}_1( \overline{\mathcal L})^2\cdot  \widehat{Z}_\rho (m) =\,{\rm ht}_{\overline{\mathcal L}}(T(m)) 
+  \int_X \widetilde{\Xi }_\rho (v,z,m)c_1( \overline{\mathcal L})^2 .\nne
\end{align}
Observe that $c_1( \overline{\mathcal L})$ is proportional to the hyperbolic measure and, as evaluated later (Remark \ref{Remvolhyp}), one has 
for the volume element
\begin{align}\label{hypvol}
c_1( \overline{\mathcal L})^2 = (18/\pi ^2)d\mu (z) = (18/\pi ^2)\frac{dx_1dy_1}{y_1^2} \frac{dx_2dy_2}{y_2^2}.\nne
\end{align}
Now from Theorem 7.61 in [BKK] p.81 we already know\\

{\bf Proposition.}\label{Propht} \emph{The Faltings height of $T(m)$ is given by}
\begin{align*}\label{eqht}
{\rm ht}_{\overline{\mathcal L}}(T(m)) =  
\begin{cases} 0 & \textrm{ if }\,\, m<0 \\
24^2 \big((\sigma (m)((1/2)\zeta (-1)+\zeta '(-1))\\\quad\quad\quad\quad + \sum_{d|m}(\frac{d \log d}{24}-\frac{\sigma (m) \log m}{48})\big) 
& \textrm{ if }\,\,m>0. \nne
\end{cases}
\end{align*}



Therefore we 
need only to study the integrals
\begin{align*}
 \int_X \widetilde{\Xi }_\rho (v,z,m)c_1(\overline{\mathcal L})^2 =  \int_X \Xi (v,z,m)c_1(\overline{\mathcal L})^2  + \int_X \rho (z)\check{\Xi } (v,z,m)c_1(\overline{\mathcal L})^2. 
\end{align*}

For the integrals $\int_X \rho (z)\check{\Xi } (v,z,m) c_1(\overline{\mathcal L})^2$ we first recall from Proposition 4.3 in Part I that 
by adding an appropriate zeroth coefficient the $q$-series
\begin{align*}
\check{\Xi }^+ (\tau ,z) = \check{\Xi }^+(v,z,0) + \sum_{m \neq 0} \check{\Xi } (v,z,m)q^m
\end{align*}
 is a modular 
form with respect to $\SL_2(\mathbb Z)$. Thus, the existence of those integrals implies our first result:\\

{\bf Theorem A.}\label{ThA}  \emph{
There exists a non-holomorphic modular form $f_\rho(\tau)$ of weight $2$ for ${\rm SL}_2(\mathbb Z)$ such that
\begin{align*}  
 f_\rho(\tau) =   \int_X \rho (z) \check{\Xi }^+ (\tau,z)c_1(\overline{\mathcal L})^2. 
  \end{align*}}

\nn {\bf Remark.} We observe that   the existence 
of the integral of $\widetilde{\Xi }_\rho $ is guaranteed by Arakelov theory. 
 Then the existence of the integral of $\int_X \rho  \check{\Xi }c_1(\overline{\mathcal L})^2 $ implies the existence of the integral 
of the Kudla Green function $\Xi $.\\

Using ${\rm O}(2,2)$-theory we are able to calculate the remaining integrals (see Theorem \ref{ThGreenint}).\\

{\bf Theorem B.}\label{ThB}  \emph{
 We have }
\begin{align*}
\int_X &\Xi (v,z,m))c_1(\overline{\mathcal L})^2 \\
& = \begin{cases}12\sigma _1(m)(1/(4\pi mv))&{\rm for} \,\,m>0\\\nonumber
 12\sigma _1(m)((1/(2\pi|m|v))e^{-4\pi |m|v} + {\rm Ei}(-4\pi |m|v)) &{\rm for}\,\,m<0.\nne
 \end{cases}
\end{align*}

It is now a pleasant exercise to relate 
these arithmetic intersection numbers for $m\neq0$ to the  
Fourier coefficients  of the Eisenstein series 
as in our Main Theorem (see Theorem \ref{ThderEiss}). 
Now, since we already know that the right hand side of our main
theorem is a non-holomorphic modular form for ${\rm Sl}_2(\mathbb Z)$ 
of weight $2$, the remaing arithmetic intersection number 
for $m=0$ must equal the coefficients of the modular form of the 
right hand side\footnote{We had spent much effort to calculate 
this identity directly, but we had not been able to do so and would be thankful for any helpful hints.}  
 \\

As already stated in Part I, our treatment of this topic owes a lot to discussions with J.Bruinier, J. Funke and S. Kudla. 
And this time we even got some local help by hints from H. Br\"uckner and J. Michali\c{c}ek. We thank them all.\\ 
 \section{Eisenstein series and its derivatives}

We take over classical material from Zagier's article \cite{Za} p.32f.
For $\tau = u+iv \in \mathbb H$ and $s \in \mathbb C$ with $\Rea s > 1$ one has the analytic Eisenstein series
\begin{align*}
E(\tau ,s) :=& (1/2) {\sum }'_{c,d}\frac{v^s}{\mid c\tau +d \mid ^{2s}}\\
=& v^s \zeta (2s) + v^s \sum_{c\in\N}\sum _{d\in\Z} \mid c\tau + d \mid ^{-2s}
\end{align*}
resp. in normalized version
\begin{align*}\label{Estar}
E^\ast (\tau ,s) := \pi ^{-s}\Gamma (s)E(\tau ,s).\nne
\end{align*}
With
$$
\zeta ^\ast (s) := \pi^{-s/2}\Gamma (s/2) \zeta (s)
$$
Zagier states the following Fourier development
\begin{align*}\label{EisFourdev}
E^\ast (\tau ,s) =& v^s \zeta ^\ast (2s) + v^{1-s} \zeta ^\ast (2s-1) \\\nonumber
&+ 2v^{1/2} \sum _{n\in\mathbb Z, n\not=0} \sigma ^\ast _{s-(1/2)} (\mid n \mid ) K_{s-(1/2)}(2\pi |n|v) e^{2\pi inu}\nne
\end{align*}
where
\begin{align*}\label{sigast}
 \sigma ^\ast _\nu (n) := \mid n \mid ^\nu \sum_{d \mid n} d^{-2\nu } =  \,\sigma ^\ast _{-\nu }(n)\nne
\end{align*}
is an entire function in $\nu $ and the {\it K-Bessel function}
$$
K_\nu (t) := \int_0^\infty e^{-t\cosh u} \cosh (\nu u) du = \,\,K_{-\nu} (t)
$$
is entire in $\nu $ and exponentially small in $t$ as $t \mapsto \infty .$

We introduce\\
\nn {\bf Definition.}
\begin{align*}\label{derEis}
E_2(\tau ,s) := (1/(2\pi i)) \partial _\tau E^\ast (\tau ,s) = (-1/(4\pi )) (\partial _v + i\partial _u) E^\ast (\tau ,s)\nne
\end{align*}
and want to study its Taylor expansion at $s=1.$ More precisely, we slightly extend the Theorem in the Introduction.\\

\nn {\bf Theorem.} \label{ThderEiss} \emph{We have
\begin{align*}
E_2(\tau ,1) = \sum_{m\in\Z}a(v,1,m) q^m\nne
\end{align*}
and, denoting by $E'_2(\tau ,s)$ the derivative of $E_2(\tau ,s)$ with respect to $s,$ we get
\begin{align*}
E'_2(\tau ,1) = \sum_{m\in\Z}a'(v,1,m) q^m\nne
\end{align*}
with}
\begin{align*}
a(v,1,m) &=
\begin{cases} \sigma (m) = \sum_{d|m} d\,\,&{\rm for} \,\,m>0\\\nonumber
- 1/24 + 1/(8\pi v)&{\rm for} \,\,m=0\\\nonumber
 0\,\,&{\rm for}\,\,m<0\nonumber
 \end{cases}\\
a'(v,1,m) &= \begin{cases} \sigma (m)(1/(4\pi mv) + \sigma '(m)/\sigma (m)) \,\,&{\rm for} \,\,m>0\\\nonumber
 -(1/24)(24\zeta '(-1)+\gamma -1+\log(4\pi v))\\\nonumber
 \quad\quad\quad\quad\quad\quad \,-(1/(8\pi v))(-\gamma +\log(4\pi v))&{\rm for} \,\,m=0\\\nonumber
 \sigma (m)(\Ei(-4\pi |m|v) + 1/(4\pi |m|v)e^{-4\pi |m|v})\,\,&{\rm for}\,\,m<0.\nne
 \end{cases}
\end{align*}


{\bf Proof.}
We have (see for instance Iwaniec \cite{Iw} p.205)
\begin{align*}
K_\nu  (t) :=& \int_0^\infty e^{-t\cosh u} \cosh (\nu u) du \\[.3cm]
=& \frac{\sqrt \pi (t/2)^\nu }{\Gamma (\nu +(1/2))} \int_1^\infty e^{-tr} (r^2-1)^{ \nu -(1/2)}dr .
\end{align*}
Hence, from(\ref{EisFourdev}) we get
\begin{align*}
E^\ast (\tau ,s) = & v^s \zeta ^\ast (2s) + v^{1-s} \zeta ^\ast (2s-1) \\[.3cm]
& + \sum _{m\in\Z, m\not=0} 2\sigma ^\ast _{s-(1/2)} (|m|) \frac{(v\mid m \mid \pi)^s}{\Gamma (s)\sqrt{|m| }}
 \int_1^\infty e^{-2\pi\mid m \mid vr} (r^2-1)^{s - 1}dr \,e^{2\pi imu}.
 \end{align*}
We abbreviate
\begin{align*}
c_0(v,s) :=& \,\,v^s \zeta ^\ast (2s) + v^{1-s} \zeta ^\ast (2s-1) ,\\[.1cm]
\check c_0(v,s) :=& \,\,\partial _vc_0(v,s) = sv^{s-1} \zeta ^\ast (2s) + (1-s)v^{-s} \zeta ^\ast (2s-1) ,\\
 \end{align*}
 and, for $m\not=0,$
\begin{align*}
c_m(v,s) :=& \,2\sigma ^\ast _{s-(1/2)} (|m|) \frac{(v|m|\pi)^s}{\Gamma (s)\sqrt{|m| }},\\
I_m (v,s):=&  \int_1^\infty e^{-2\pi\mid m \mid vr} (r^2-1)^{s - 1} dr ,\\[.3cm]
J_m (v,s):=&  \int_1^\infty e^{-2\pi\mid m \mid vr} (r^2-1)^{s - 1}r dr ,\\[.1cm]
 \end{align*}
 and get
 \begin{align*}
 E_2(\tau ,s) = - (1/(4\pi))( \check c_0(v,s) + \sum '_m ((s/v \,\,-& 2\pi m)c_m(v,s)I_m(v,s) \\[.1cm]
-& 2\pi |m| c_m(v,s)J_m(v,s))e(mu)
 \end{align*}
  i.e.
 \begin{align*}
 E_2(\tau ,s) =& - \frac{1}{4\pi }(\check c_0(v,s)\\\nonumber
&+ \sum_{m>0} ((s/v)c_m(v,s)I_m(v,s) - 2\pi |m| c_m(v,s)(I_m(v,s)+J_m(v,s))e(mu)\\\nonumber
&+ \sum_{m<0} ((s/v)c_m(v,s)I_m(v,s) + 2\pi |m| c_m(v,s)(I_m(v,s)-J_m(v,s))e(mu)\,).\nne
 \end{align*}
Using Maple, one can determine from here the first two terms of the Taylor expansion of each coefficient of $e(mu)$ 
and hence get the claims in the Theorem. For those who don't like Maple, we give a direct proof in an appendix.\\
 
\hfill$\Box$\\

\nn {\bf Remark.}\label{remsig} The sigmas in this calculations are those from the paper by Zagier
$$
\sigma ^\ast _s (n) := |n|^s \sum_{d \mid n, d>0} d^{-2s } =  \,\sigma ^\ast _{-s }(n).
$$
Hence one has
\begin{align*}
\sqrt m \,\sigma ^\ast  _{1/2}(m) = \sum_{d|m} d = \sigma (m). 
\end{align*}
We set
\begin{align*}\label{eqsig}
\sigma '(m)/\sigma (m) &:= {\sigma ^\ast }'(m)_{1/2}/\sigma ^\ast  _{1/2}(m) \\\nne
&=  (\sigma (m)\log m -2 \sum_{d|m}d\log d)/\sigma (m).\\\nonumber
\end{align*}
As an immediate consequence to our 
 Theorem, for the coefficients \eqref{eqmodEis} of the modified Eisenstein series $\mathbb E'_2(\tau ,1) ,$ we get\\
\nn{\bf Corollary.} \label{CormodA} One has
\begin{align*}
A'(v,1,m) &= -12\begin{cases} \sigma (m)(4(\zeta '(-1)/\zeta (-1)+1/2)-1/(4\pi mv)\\\nonumber
\quad \quad\quad \quad\quad \quad\quad \quad\quad \quad + {\sigma ^\ast }'_{1/2}(m)/\sigma ^\ast _{1/2}(m))) \,\,&{\rm for} \,\,m>0\\\nonumber
 \,3\zeta '(-1)-(1/8)+(\gamma /24) + (1/24)\log(4\pi v)\\\nonumber
\quad + (1/8\pi v)(-48\zeta '(-1)-\gamma +2+\log(4\pi v)). &{\rm for}\,\,m=0 \\\nonumber
 \sigma (|m|)(\Ei(-4\pi |m|v) + 1/(4\pi |m|v)e^{-4\pi |m|v})\,\,&{\rm for}\,\,m<0\nne
 \end{cases}
\end{align*}

\section{Boundary function integral}
In Section 2 in Part I we introduced a partition of the unity $\rho $ with respect to the boundary $D$ and 
the boundary function $\check{\Xi }^+ (\tau ,z)$ 

\begin{align*}
\check{\Xi }^+ (\tau ,z) &= \sum_m \check \Xi (v,z,m)q^m - (1/2v)t(s+1/s)\\\nonumber
\check \Xi (v,z,m) &= (1/2)\sum_{-bc=m} \check \xi (v,z;b,c)\nonumber
\end{align*}
with $t=\sqrt{y_1y_2}, s=\sqrt{y_1/y_2}$ (unfortunately we here have the same letter as the one denoting the 
variable in the zeta and Eisenstein series but the kind reader will know to make the difference) and
\begin{align*}
\check \xi (v,z,;b,c) &=  \,(t/\sqrt v)\,\, ( \operatorname{B}(v,s;b,c) - \operatorname{I}(v,s;b,c)) \\ \nonumber
\operatorname{B}(v,s;b,c) &=\int_1^\infty e^{-\pi v(b/s+cs)^2r}r^{-3/2}dr\\\nonumber
\operatorname{I}(v,s;b,c) &= \begin{cases}  4 \pi \sqrt v\,\min (\mid bs^{-1} \mid ,\mid cs \mid )) \,\, &{\rm if}\,\,-bc > 0\\
  0  \,\,\quad \quad &{\rm if}\,\,-bc \leq  0.
  \end{cases}
\end{align*}

And in Proposition 4.3 of Part I we proved the modularity of $\check{\Xi }^+ (\tau ,z)$ as a function in $\tau .$ 
From there we come to the following result:

\nn {\bf Theorem.} \label{thm:frho}\emph {
There exists a modular form $f_\rho$ such that}
\begin{align*} \label{eq:frho}
 f_\rho (\tau ) =   \int_X \rho(z)\check{\Xi }^+ (\tau ,z))d\mu.\nne
 \end{align*}

{\bf Proof.} As we have modularity in $\tau $ of $\check{\Xi }^+$ and 
since the integral does not affect the $\tau -$variable, the
 modularity follows as soon as we checked the existence of the integrals. 
For $m=0$ this is evaluated in the Proposition \ref{Propxi0} below 
and for $m\not=0$  that follows from the three lemmata below as in these we have integrals of type
$$
\int_X  t F(s)d\mu = \int_X tF(s) dx_1dx_2dsdt/st^3
$$  
and, as the integrand does not depend on $x_1,x_2,$ once the $s-$integration is done, 
one has a finite value as 
$$
\int_{t>t_0}dt/t^2 < \infty .
$$
\hfill $\Box$

\nn{\bf Lemma.} \emph {For $m=-bc<0$ and 
\begin{align*}
\tau (m) := \sharp \{d>0:d \,| \,|m|\}
\end{align*}
one has}
\begin{align*}\sum_{-bc=m}\int_0^\infty B(v,s;b,c)ds/s
&\leq 2\tau (m)(1/(2\sqrt {|m|v})e^{-4\pi |m|v} + 2\pi \sqrt {|m|v}\,{\Ei}(-4\pi |m|v))
\end{align*}
{\bf Proof.} Replacing $s$ by $s\sqrt {c/b}$ we get
\begin{align*}
\int_0^\infty B(v,s;b,c)ds/s &= \int_0^\infty \int_1^\infty e^{-\pi v(b/s+cs)^2r}dr/r^{3/2}ds/s\\
&= \int_0^\infty \int_1^\infty e^{-\pi v|m|((1/s)^2+s^2) + 2)r}dr/r^{3/2}ds/s\\
\end{align*}
and with $s=e^\varphi $ and $\cosh \varphi = 1 + \varphi ^2/2 + \dots$ we estimate
\begin{align*}
\int_0^\infty B(v,s;b,c)ds/s &= \int_1^\infty \int_{-\infty }^\infty e^{-2\pi vr|m|\cosh 2\varphi }d\varphi e^{-2\pi |m|vr}dr/r^{3/2}\\
&\leq  \int_1^\infty \int_{-\infty }^\infty e^{-4\pi v|m|r\varphi ^2}d\varphi e^{-4\pi |m|vr}dr/r^{3/2}\\
&=  1/(2\sqrt {|m|v})\int_1^\infty  e^{-4\pi |m|vr}dr/r^{2} \\
&=  1/(2\sqrt {|m|v})e^{-4\pi |m|v} + 2\pi \sqrt {|m|v}\,{\Ei}(-4\pi |m|v)).
\end{align*}
\hfill$\Box$\\

\nn{\bf Lemma.} \label{LeB_m}\emph {For $m=-bc>0$ one has }
\begin{align*}
 \sum_{-bc=m}\int_0^\infty B(v,s;b,c)ds/s
\leq 2\tau (m)(1/(2\sqrt {|m|v}))
\end{align*}
{\bf Proof.} Replacing $b$ by $-b$ one has $\sum_{-bc=m} = 2\sum_{b,c>0,bc=m}$ and again  $s$ by $s\sqrt {c/b}$ we get this time
\begin{align*}
\int_0^\infty B(v,s;b,c)ds/s &= \int_0^\infty \int_1^\infty e^{-\pi v(b/s+cs)^2r}dr/r^{3/2}ds/s\\
&= \int_0^\infty \int_1^\infty e^{-\pi vm((1/s)^2+s^2) - 2)r}dr/r^{3/2}ds/s\\
\end{align*}
and with $s=e^\varphi $ 
\begin{align*}
\int_0^\infty B(v,s;b,c)ds/s &= \int_1^\infty \int_{-\infty }^\infty e^{-2\pi vrm\cosh 2\varphi }d\varphi e^{2\pi mvr}dr/r^{3/2}\\
&\leq  \int_1^\infty \int_{-\infty }^\infty e^{-4\pi v|m|r\varphi ^2}d\varphi dr/r^{3/2}\\
&=  1/(2\sqrt {mv})\int_1^\infty  dr/r^{2} \\
&=  1/(2\sqrt {mv}).
\end{align*}\hfill$\Box$\\

\nn{\bf Lemma.}\label{LeI_m} \emph {For $m=-bc>0$ one has }
\begin{align*}
 \sum_{-bc=m}\int_0^\infty \min (|b/s|,|cs|)ds/s = 4\tau (m) \sqrt m
\end{align*}
{\bf Proof.} Replacing $s$ by $s\sqrt {|c|/|b|}$ we get 
\begin{align*}
\sum_{-bc=m}\int_0^\infty \min (|b/s|,|cs|)ds/s &= \sum_{-bc=m}\sqrt m \int_0^\infty \min (1/s,s)ds/s\\
&= 2\sqrt m\tau (m)(\int_0^1 sds/s + \int_1^\infty ds/s^2)\\
&= 2\sqrt m\tau (m)\cdot 2.
\end{align*}\hfill$\Box$\\

\nn {\bf Proposition.} \label{Propxi0} \emph { 
For $m=bc=0$  one has }
\begin{align*}
 \int_X \check \Xi^+(v,z,0) d\mu < \infty
\end{align*}


{\bf Proof.} From the Remark 2.23 from Part I we get
 \begin{align*}
2\cdot \check{\Xi}(v,z,0) = &(t/\sqrt v)( \sum_{b\not=0}B(v,s;b,0) + \sum_{c\not=0}B(v,s;0,c) + B(v,s;0,0))\nonumber\\
= & -2t/\sqrt v + t(s+1/s)(1/v + (2/\pi )\zeta (2)) \nonumber\\
&-(2t/\pi)((1/s)\sum _{b\in \mathbb N}e^{-\pi s^2b^2/v}/b^2 + s\sum_{c\in \mathbb N}e^{-\pi c^2/(s^2v)}/c^2 )\nne 
\end{align*}
and from (4.2.6) of Part I
 \begin{align*}
2\cdot \check{\Xi}^+(v,z,0) = & 2\cdot \check{\Xi}(v,z,0) - (1/v)t(s+1/s)\nonumber\\
= & -2t/\sqrt v + t(s+1/s)(2/\pi )\zeta (2)) \nonumber\\
&-(2t/\pi)((1/s)\sum _{b\in \mathbb N}e^{-\pi s^2b^2/v}/b^2 + s\sum_{c\in \mathbb N}e^{-\pi c^2/(s^2v)}/c^2 ).\nne 
\end{align*}
{\bf Step 1.} We start by integrating
\begin{align*}
I' &= \int_{K_1<y_1} \int_{K_2<y_2< T} t dy_1dy_2 / (y_2^2 y_1^2)\\
&= \int_{K_2}^{T}\int_{K_1}^\infty dy_1/y_1^{3/2}dy_2/y_2^{3/2} = (4/\sqrt K_1)(1/\sqrt K_2 - 1/\sqrt T).
\end{align*}
and
\begin{align*}
I_0 &= \int_{K_1<y_1} \int_{K_2<y_2< T} y_2 dy_1dy_2 / (y_2^2 y_1^2)\\
&= \int_{K_2}^{T}\int_{K_1}^\infty dy_1/y_1^{2}dy_2/y_2\\
&=\int_{K_2}^{T} (1/K_1)dy_2/dy_2 = (1/K_1)(\log T - \log K_2).
\end{align*}
{\bf Step 2.} For $b\not=0$ we look at 
\begin{align*}
I_b :=& (1/b^2)\int_{K_2}^{T}y_2\int_{K_1}^{\infty }e^{-(\pi /v)b^2y_1/y_2}dy_1dy_2/(y_1y_2)^2\\
=& (1/b^2)\int_{K_2}^{T} ([e^{-(\pi /v)b^2y_1/y_2}(-1/y_1)]_{K_1}^{\infty }\\ 
&- \int_{K_1}^{\infty }(\pi b^2/(y_2v))e^{-(\pi /v)b^2y_1/y_2}dy_1/y_1 )dy_2/y_2\\
=& (1/b^2)\int_{K_2}^{T}([e^{-(\pi /v)b^2K_1/y_2}(1/K_1)] - \int_1^\infty (\pi b^2/(y_2v))e^{-(\pi /v)b^2K_1y_1/y_2}dy_1/y_1)dy_2/y_2.
\end{align*}
We remind
\begin{align*}
-{\rm Ei}\,(-x) = \int_{1 }^\infty  e^{-xt}dt/t = - \gamma  - \log |x| + x - x^2/(2\cdot 2!) +  \dots 
\end{align*}
and have for $T \longrightarrow \infty $ in the first term
\begin{align*}
I_{b,1} &=(1/(b^2K_1))\int_{K_2}^{T}e^{-(\pi /v)b^2K_1/y_2}dy_2/y_2 \\
&=  (1/(b^2K_1))\int_{1/T}^{1/K_2}e^{-(\pi /v)b^2K_1u}du/u\\
&=  (1/(b^2K_1))\int_{1}^{T/K_2}e^{-(\pi /v)b^2K_1u/T}du/u\\
&\leq    (1/(b^2K_1))\int_{1}^{\infty }e^{-(\pi /v)b^2K_1u/T}du/u\\
&=  (1/(b^2K_1))(- \gamma  - \log x + x +\dots)
\end{align*}
where $x = \pi (b^2/v)(K_1/T).$ 
 And for the second term in $I_b$
\begin{align*}
I_{b,2} :&= (\pi /v)\int_{K_2}^T(\int_1^\infty e^{-(\pi /v)b^2K_1y_1/y_2}dy_1/y_1)dy_2/y_2^2\\
&=-(\pi /v)\int_{K_2}^T \Ei(-(\pi /v)b^2K_1/y_2)dy_2/y_2^2\\ 
&=-(\pi /v)\int_{1/T}^{1/K_2} \Ei(-(\pi /v)b^2K_1u)du
\end{align*}
with $\alpha = \pi b^2K_1/v$ we have
\begin{align*}
I_{b,2} =& (\pi /v)\int_{K_2}^T \int_1^\infty e^{-\alpha y_1/y_2}dy_1/y_1dy_2/y_2^2\\
=& (\pi /v)\int_1^\infty \int_{1/T}^{1/K_2}[e^{-\alpha y_1u}dudy_1/y_1\\
=& (\pi /(v\alpha ))\int_1^\infty [e^{-\alpha y_1/T} - e^{-\alpha y_1/K_2}]dy_1/y_1^2\\
=& (1/(b^2K_1))\int_1^\infty [e^{-\alpha y_1/T} - e^{-\alpha y_1/K_2}]dy_1/y_1^2\\
=& (1/(b^2K_1))([e^{-\alpha /T} - e^{-\alpha /K_2}])\\ 
&- (\pi /v) \int_1^\infty [e^{-\alpha y_1/T}/T - e^{-\alpha y_1/K_2}/K_2]dy_1/y_1\\
=& (1/(b^2K_1))([e^{-\alpha /T} - e^{-\alpha /K_2}]) \\
&- (\pi /(vK_2))\Ei(-\pi (b^2/v)K_1/K_2) + (\pi /(vT))\Ei(-\pi (b^2/v)K_1/T)
\end{align*}
i.e. something finite for $T \longrightarrow \infty $ as the first terms are harmless and for the last one one has
\begin{align*}
(\pi /(vT))\Ei(-\pi (b^2/v)K_1/T) =&  (\pi /(vT))(- \gamma  - \log (-\pi (b^2/v)K_1/T)\\
&+ (-\pi (b^2/v)K_1/T) +\dots.
\end{align*}
with $(1/T)\log T \rightarrow 0.$\\

{\bf Step 3.} We remark that for $T \longrightarrow \infty $ $I_0/b^2$  and $I_b$ have the same singularity, namely 
$1/(b^2K_1)\log T.$\\

{\bf Step 4.} The same way, we have the the same singularity coming from
\begin{align*}
I'_0 &= \int_{K_2<y_2} \int_{K_1<y_1< T} y_1 dy_1dy_2 / (y_2^2 y_1^2)\\
&= \int_{K_1}^{T}\int_{K_2}^\infty dy_2/y_2^{2}dy_1/y_1\\
&=\int_{K_1}^{T} (1/K_2)dy_1/dy_1 = (1/K_2)(\log T - \log K_2).
\end{align*}
and
\begin{align*}
I_c :=& (1/c^2)\int_{K_1}^{T}y_1\int_{K_2}^{\infty }e^{-(\pi /v)c^2y_2/y_1}dy_1dy_2/(y_1y_2)^2\\
=& (1/c^2)\int_{K_1}^{T} ([e^{-(\pi /v)c^2y_2/y_1}(-1/y_1)]_{K_2}^{\infty }\\ 
&- \int_{K_2}^{\infty }(\pi c^2/(y_1v))e^{-(\pi /v)c^2y_2/y_1}dy_2/y_2 )dy_1/y_1\\
=& (1/c^2)\int_{K_1}^{T}([e^{-(\pi /v)c^2K_2/y_1}(1/K_2)]\\
 &- \int_1^\infty (\pi c^2/(y_1v))e^{-(\pi /v)c^2K_2y_2/y_1}dy_2/y_2)dy_1/y_1.
\end{align*}
Hence all together adds up to something finite.\hfill$\Box$\\

\section{Kudla's Green function integral for $m\neq 0$}

At first we remark that, as explained at the end of the Introduction, it follows immediately from Theorem \ref{thm:frho} 
that for $m\ne 0$ the integrals in question exist.\\

We look at the Green function integral
\nn {\bf Definition.}
\begin{align*}\label{Greenint}
I_m :=& \int_{\Gamma \backslash \Ha \times \Gamma \backslash \Ha} \Xi (v,z,m) d\mu (z)\\\nonumber
=& \int_{\Gamma \backslash \Ha \times \Gamma \backslash \Ha} (1/2)\sum_{M\in L_m}\xi (v,z,m) d\mu (z)\nne
\end{align*}

and want to prove the following\\

\nn {\bf Theorem.} \label{ThGreenint} \emph {With $\sigma (m) := \sum _{d|m}d$ one has }
\begin{align*}\label{eqGreenint}
I_m  =\begin{cases}  \sigma (|m|)(\pi /3)(1/(2v|m|))(e^{-4\pi v|m|} + 4\pi v|m|{\Ei}(-4\pi |m|v))\,\,&{\rm for}\,\,m<0\\
 \sigma (m)(\pi /3) (1/(2vm)) \,\,&{\rm for}\,\,m>0.\end{cases}\nne
\end{align*}

{\bf Proof.} At first we assemble some tools. For $m \in \mathbb N$ and 
\begin{align*}
\Gamma &= \SL(2,\mathbb Z),\\
L_m &= \{M = \begin{pmatrix}a&b \\c&d \end{pmatrix} \in M_2(\mathbb Z);\, \det M = m \},\\
L_m^\ast &= \{ M \in L_m; \,M \,\,{\rm primitive} \},
\end{align*}
one has the standard facts (see for instance Ogg's book \cite{Og} p.II-7 and IV-4)
\begin{align*}\label{eqLm}
L_m &= \cup _{ad=m,d>0,a\mid d}\,\Gamma \begin{pmatrix}a& \\&d \end{pmatrix} \Gamma \\\nonumber
&= \cup _{ad=m,d>\,0,\,b\!\!\mod d}\,\Gamma \begin{pmatrix}a& b\\&d \end{pmatrix} \nne
\end{align*}
and
\begin{align*}\label{eqLastm}
L^\ast _m &= \Gamma  \begin{pmatrix}m& \\&1 \end{pmatrix}\Gamma = \,\,\sqcup_\alpha  \Gamma \alpha \nne
\end{align*}
with
$$
\alpha = \begin{pmatrix}a&b \\&d \end{pmatrix},\,\, ad = m, d > 0, 0 \leq b < d, (a,b,d) = 1.
$$
One has
$$
[L^\ast _m : \Gamma ] = m \prod _{p\mid m} (1+(1/p)) =: \psi (m)
$$
and also
$$
[\Gamma :\Gamma _0(m)] = \psi (m).
$$
Moreover, we have
$$
[L_m : \Gamma ] = \sum _{d\mid m} d = \sigma (m)
$$
and hence the formula 
(which is easily verified using the multiplicativity of $\sigma $, 
see for instance Rankin \cite{Ra} p. 285)
\begin{align*}\label{eqRan}
\sigma (m) = \sum _{n^2\mid m} \psi (m/n^2).\nne
\end{align*}
$\Gamma $ acts transitively by right multiplication 
on $\Gamma \backslash L^\ast _m$ with isotropy group $\Gamma _0(m)$ at the coset 
$\Gamma \begin{pmatrix}m& \\&1 \end{pmatrix}$ (for example as in Knapp \cite{Kn} p.256, Proposition 9.3). Hence, one has a bijection
\begin{align*}\label{eqLmbij}
L_m^\ast \simeq  \Gamma \begin{pmatrix}m& \\&1 \end{pmatrix}(\Gamma _0(m)\setminus \Gamma ).\nne
\end{align*}
It also is a standard fact that one has
\begin{align*}\label{eqvol}
{\rm vol }(\Gamma \backslash \mathbb H) = \int_{\Gamma \backslash \mathbb H } \frac{dxdy}{y^2} =\pi/3.\nne
\end{align*}
and for $m>0$ (e.g. \cite{FB} p.375)
\begin{align*}\label{eqvolm}
{\rm vol }(\Gamma _0(m)\backslash \mathbb H) = \psi (m)\pi/3.\nne
\end{align*}


After the preparation of these tools, we come to calculate the Green function integral
$$
I_m := \int_{\Gamma \backslash \mathbb H \times \Gamma \backslash \mathbb H} \Xi (v,z,m) d\mu (z)
$$
with
$$
d\mu (z) = d\mu (z_1)d\mu (z_2) = \prod _{j=1,2} \frac{dx_jdy_j}{y_j^2}.
$$
We do this in several steps.\\

{\bf Step 1. The integral for squarefree positive $m$}\\

To simplify things, we start by treating the special case of squarefree  $m > 0$ where $L^\ast _m = L_m.$ \\

As $\overline \Gamma = \Gamma \times \Gamma $ acts on $L_m$ via $M \mapsto \gamma _1M{}^t\gamma _2 =: M^\gamma ,$ we have
\begin{align*}
I_m =& \int_{\Gamma \backslash \mathbb H \times \Gamma \backslash \mathbb H} (1/2)\sum_{M \in L^\ast _m} \xi (v,z,M) d\mu (z)\\[.1cm]
=&\int_{\Gamma \backslash \mathbb H \times \Gamma \backslash \mathbb H}(1/2)
\sum_{\gamma _1  \in \Gamma , \beta \in (\Gamma _0(m)\setminus \Gamma )} 
\xi (v,(z_1, z_2);\gamma _1\begin{pmatrix}m& \\&1 \end{pmatrix}{}^t\beta  ) d\mu (z)\\[.1cm]
=&\int_{\Gamma \backslash \mathbb H \times \Gamma \backslash \mathbb H}(1/2)
\sum_{\gamma _1  \in \Gamma , \beta \in (\Gamma _0(m)\setminus \Gamma )} 
\xi (v,(\gamma _1 ^{-1}z_1, \beta ^{-1}z_2);\begin{pmatrix}m& \\&1 \end{pmatrix} ) d\mu (z)
\end{align*}
where we use the homogenity $\xi (gz,M^g) = \xi (z,M).$ Hence, one has
$$
I_m = (1/2)\int_{ \mathbb H \times (\Gamma _0(m)\backslash \mathbb H)}  
\xi (v,z,\begin{pmatrix}m& \\&1 \end{pmatrix}) d\mu (z)
$$
with
$$
\xi (v,z,\begin{pmatrix}m& \\&1 \end{pmatrix}) =  \int _1^\infty e^{-2\pi v
R(z,(\begin{smallmatrix}m& \\&1 \end{smallmatrix}))u}du/u
$$
and
$$
R(z,\begin{pmatrix}m& \\&1 \end{pmatrix} = (1/(2y_1y_2))|m + z_1z_2|^2.
$$
We simplify this by changing two times our coordinates. \\

i) The change
$$
(z_1,z_2) \longmapsto (z_1,mz_2)
$$
leads to $d\mu (z) \mapsto d\mu (z)$ and 
$$
R(z,( \begin{smallmatrix}m& \\&1 \end{smallmatrix})) \longmapsto m R(z,( \begin{smallmatrix}1& \\&1 \end{smallmatrix})).
$$
ii) For $g_{z_2}$ with $g_{z_2}(i) = z_2$ we take
$$
z = (z_1,z_2) \longmapsto g(z) = ({}^tg_{z_2}(z_1) ,g_{z_2}^{-1}(z_2)) =: ({z_1}', i)
$$
and have
\begin{align*}
R(z,( \begin{smallmatrix}1& \\&1 \end{smallmatrix})) =& R(g(z),( \begin{smallmatrix}1& \\&1 \end{smallmatrix})^g)\\[.1cm]
=& R({z_1}', i;,( \begin{smallmatrix}1& \\&1 \end{smallmatrix}))\\[.1cm]
=& (1/(2{y_1}'))\mid 1 + i{z_1}' \mid ^2
\end{align*}
and finally
$$
\begin{array}{rl}
I_m =& {\rm vol}(\Gamma _0(m) \setminus \mathbb H)(1/2) \int_\mathbb H (\int_1^\infty e^{-2\pi vm R({z_1}',i ; ( \begin{smallmatrix}1& \\&1 \end{smallmatrix}))u}du/u) d\mu ({z_1}')\\[.3cm]
=& {\rm vol}(\Gamma _0(m) \setminus \mathbb H) (1/2)\int_\mathbb H (\int_1^\infty e^{-\pi vm(1/(y_1))((1-y_1)^2 + x_1^2) u}du/u)d\mu (z_1).\\
\end{array}
$$
Thus one is reduced to a two-dimensional integral
\begin{align*}\label{eqI'm}
I'_m = \int_\mathbb H\int_1^\infty e^{-\pi vm(1/(y_1))((1-y_1)^2 + x_1^2) u}du/u)d\mu (z_1).\nne
\end{align*}
Using parts of 
the $\SO(1,2)-$theory as for instance in Brunier-Funke \cite{BF}, we change coordinates
\begin{align*}
\mathbb R^2 \longrightarrow \mathbb H,\,\,(r,\varphi )\longmapsto (z=x+iy)
\end{align*}
with
$$
y = 1/(\cosh r - \sinh r \cos \varphi ) , \,\, x = - \sinh r \sin \varphi /(\cosh r - \sinh r \cos \varphi ) .
$$
i.e.,
\begin{align*}
(x^2+y^2+1)/(2y) &= \cosh r\\
 (x^2+y^2-1)/(2y) &= \sinh r \cos \varphi \\
 -x/y &= \sinh r \sin \varphi .
\end{align*}
A small calculation shows that one has
$$
d\mu (z) = dx \wedge dy/y^2 = \sinh r \,\,dr \wedge d\varphi .
$$
We get
\begin{align*}
I'_m  &=
\int_\mathbb H (\int_1^\infty e^{-\pi vm(1/y)((1+x^2 + y^2 - 2y)u}du/u)d\mu (z)\\[.1cm]
&= \int _0^{2\pi} \int_0^\infty (\int_1^\infty e^{-2\pi vm(\cosh r - 1)u}du/u)\sinh r dr d\varphi \\[.1cm]
&= 2\pi \int_1^\infty \int_1^\infty e^{-2\pi vmu(t-1)}dt du/u\\[.1cm]
&= 2\pi \int_1^\infty (\int_1^\infty e^{-2\pi vmut}dt)e^{\pi vmu} du/u\\[.1cm]
&= 2\pi \int_1^\infty [e^{-2\pi vmut}/(-2\pi vmu)]^\infty _1\,e^{2\pi vmu} du/u\\[.1cm]
&= (1/vm)\int_1^\infty u^{-2}du\\[.2cm]
&= (1/vm) .
\end{align*}
and hence
\begin{align}
I_m = {\rm vol}(\Gamma _0(m) \setminus \mathbb H) (1/2)(1/vm) = \sigma (m)\pi /(6vm).\nne
\end{align}

{\bf Step 2. The integral for squarefree negative $m$}\\

For negative $m$ we need slight changes in the calculation of the integral $I_m$ in the second part of our proof. 
At first one has
$$
L_m^\ast \simeq  \Gamma \begin{pmatrix}m& \\&1 \end{pmatrix}(\Gamma _0(|m|)\setminus \Gamma ).
$$
Then the transformation
$$
(z_1,z_2) \longmapsto (z_1,mz_2)
$$
transforms $\mathbb H \times \mathbb H$ to $\mathbb H \times \overline {\mathbb H}.$ Hence, in the next step, we have to replace 
the old $g_{z_2}$ by another one with $g_{z_2}(-i) = z_2$ and come to 
\begin{align*}\label{eqImneg}
I_m =& {\rm vol}(\Gamma _0(|m| ) \setminus \mathbb H) (1/2)\int_\mathbb H (\int_1^\infty e^{-2\pi v\mid m \mid R({z_1}',-i ; ( \begin{smallmatrix}1& \\&1 \end{smallmatrix}))u)}du/u) d\mu ({z_1}')\\[.1cm]\nonumber
=& {\rm vol}(\Gamma _0(|m| ) \setminus \mathbb H) \,(1/2)
\int_\mathbb H (\int_1^\infty e^{-\pi(1/y_1)((1+y_1)^2 + x_1^2)v\mid m \mid u}du/u)d\mu (z_1)\\[.1cm]\nonumber
=& {\rm vol}(\Gamma _0(|m| ) \setminus \mathbb H) \,(1/2)
\int_\mathbb H (\int_1^\infty e^{-\pi(1/y_1)((1+x_1^2+y_1^2 + 2y_1)v\mid m \mid u}du/u)d\mu (z_1)\\[.1cm]\nonumber
=& {\rm vol}(\Gamma _0(|m| ) \setminus \mathbb H) \,(1/2)
2\pi \int_0^\infty  (\int_1^\infty e^{-2\pi v\mid m \mid (\cosh r + 1)u}du/u)\sinh r dr\\[.1cm]\nonumber
=& {\rm vol}(\Gamma _0(|m| ) \setminus \mathbb H) \,(1/2)
2\pi \int_1^\infty  (\int_1^\infty e^{-2\pi v\mid m \mid (t + 1)u}dt)du/u\\[.1cm]\nonumber
=& {\rm vol}(\Gamma _0(|m| ) \setminus \mathbb H) \,(1/2)
(-1/vm )\int_1^\infty e^{-4\pi v\mid m \mid u}u^{-2}du\\[.1cm]\nonumber
=& {\rm vol}(\Gamma _0(|m| ) \setminus \mathbb H) \,(1/2)
((-1/vm) e^{-4\pi v\mid m \mid } - 4\pi \int_1^\infty e^{-4\pi v\mid m \mid u}du/u).\nne
\end{align*}

{\bf Step 3. The integral for general $m\not=0$}\\

We use the results which we already have and for positive $m$ calculate
\begin{align*}
I_m =& \int_{\Gamma \backslash \mathbb H \times \Gamma \backslash \mathbb H} (1/2)\sum_{M \in L_m} \xi (z,M) d\mu (z)\\[.1cm]
 =& \int_{\Gamma \backslash \mathbb H \times \Gamma \backslash \mathbb H} \sum_{n^2\mid m}(1/2)\sum_{M \in L^\ast _{m/n^2}} \xi (z,nM) d\mu (z)\\[.1cm]
 =& \sum_{n^2\mid m}\int_{\Gamma \backslash \mathbb H \times \Gamma \backslash \mathbb H} (1/2)\sum_{M \in  L^\ast _{m/n^2}}  \xi (z,nM) d\mu (z).
\end{align*}

With
$$
 L^\ast _{m/n^2} = \Gamma  \begin{pmatrix}m/n^2& \\&1 \end{pmatrix}(\Gamma _0(m/n^2)\setminus \Gamma )
$$
we get
$$
\begin{array}{rl}
I_m =& (1/2)\sum_{n^2\mid m}\int_{\mathbb H}\int_{\Gamma _0(m)\backslash \mathbb H}  
\xi (z_1, z_2;n\begin{pmatrix}m/n^2& \\&1 \end{pmatrix}) d\mu (z).
\end{array}
$$
One has
$$
R(z,n\begin{pmatrix}m/n^2& \\&1 \end{pmatrix} ) = (1/2y_1y_2)\mid m/n + nz_1z_2\mid^2
$$
and changing the coordinates, as in the second part of our proof above, this time by $(z_1,z_2) \longmapsto (z_1,mz_2/n^2)$ we get 
$$
R(z,n\begin{pmatrix}m/n^2& \\&1 \end{pmatrix} ) \longmapsto (1/2y_1y_2)m\mid 1+z_1z_2 \mid^2 = mR(z,\begin{pmatrix}1& \\&1 \end{pmatrix})
$$  
and hence, as in the second part above, via $z \longmapsto g(z) =(z_1,i)$
$$
\begin{array}{rl}
I_m =& (1/2)\sum_{n^2\mid m}{\rm vol}(\Gamma _0(|m/n^2|) \setminus \mathbb H)
\int_{\mathbb H}\int_1^\infty e^{-2\pi mR(z_1,i;(\begin{smallmatrix}1& \\&1 \end{smallmatrix}))u}(du/u) d\mu (z_1)\\[.3cm]
=& (1/2)\sum_{n^2\mid m}{\rm vol}(\Gamma _0(|m/n^2|) \setminus \mathbb H) I'_m \\[.3cm]
=& (1/2)\sum_{n^2\mid m}\psi (m/n^2) I'_m = \sigma _1(m)\kappa \,I'_m \\[.3cm]
=& \sigma _1(m)(\pi /3)(1/2)(1/vm)
\end{array}
$$
where we used the formulae (\ref{eqLm}) and (\ref{eqRan}) from the first part of the proof
 and put $\kappa := {\rm vol}(\Gamma \setminus \mathbb H) = \pi/3.$ \\

 For negative $m$ one gets the same way,
 with $m$ replaced by $|m|,$ analogously the formula (\ref{eqImneg}) from above
\begin{align*}
I_m = \sigma _1(|m|)(\pi /3)(1/v|m|)(e^{-4\pi v|m|} + 4\pi v|m|{\Ei}(-4\pi |m|v)).
\end{align*}\hfill$\Box$\\



\section{Proof of the Main Theorem}

Now we  relate the results obtained for $m\not=0$ in Theorem \ref{ThGreenint} for the 
Green function integral $I_m$ to the Fourier coefficients of our modified Eisenstein series in Corollary \ref{CormodA}.\\

\nn {\bf Theorem.} \label{CorKucon1} \emph {For $m\not=0,$ with }(\ref{cxcxZm}) \emph {one has}
\begin{align*} \label{eqKucon1}
\widehat{c}_1( \overline{\mathcal L})^2\cdot  \widehat{Z} (m) &= {\rm ht}_{\overline{\mathcal L}}(T(m))
+  \int_X \widetilde{\Xi }(v,z,m))c_1(\bar{\mathcal L})^2  \\\nonumber
&= \,A'(v,1,m).\nne
\end{align*}

{\bf Proof.} We already observed
\begin{align*} 
 A'(v,1,m) = -48(\frac{\zeta '(-1)}{\zeta (-1)} + \frac{1}{2})a(v,1,m) + 12 a'(v,1,m).
\end{align*} 
Hence, for $m > 0,$ from Theorem \ref{ThderEiss} resp. Corollary \ref{CormodA} and the Remark \ref{remsig} on the different sigmas one has
\begin{align*} 
 A'(v,1,m) &= -48\sigma (m)(\frac{\zeta '(-1)}{\zeta (-1)} + \frac{1}{2})\\
&\quad\quad + 12\sigma (m)( 1/(4\pi mv) 
+ {\sigma ^\ast }'(m)_{1/2}/\sigma ^\ast  _{1/2}(m))\\
&=   -48\sigma (m)(\frac{\zeta '(-1)}{\zeta (-1)} + \frac{1}{2})\\ &\quad \quad + 12(\sigma (m)/(4\pi mv)+(\sigma (m)\log m -2 \sum_{d|m}d\log d))  .
\end{align*}
On the other hand, as already remarked in the introduction, from Theorem 7.62 in \cite{BKK} and Theorem \ref{ThGreenint}

\begin{align*}
\widehat{c}_1( \overline{\mathcal L})^2\cdot  \widehat{Z} (m) =& \,\,{\rm ht}_{\overline{\mathcal L}}(T(m)) + c \int_X \widetilde{\Xi }(v,z,m))d\mu \\\nonumber
=& \,\, 24^2 (\sigma (m)((1/2)\zeta (-1)+\zeta '(-1))+ \sum_{d|m}(\frac{d \log d}{24}-\frac{\sigma (m) \log m}{48}))\\
&+ c \sigma _1(m)(\pi /3)2\pi (1/(4\pi mv)).
\end{align*}


For $m < 0$ one has by Theorem \ref{ThderEiss} $a(v,1,m) = 0$ and
\begin{align*} 
 A'(v,1,m) = -(\Ei(-4\pi |m|v)+e^{-4\pi |m|v}/(4\pi |m|v))\sigma _1(m))
\end{align*} 
and again from Theorem \ref{ThGreenint}
\begin{align*} 
\widehat{c}_1( \overline{\mathcal L})^2\cdot  \widehat{Z} (m) &= \,{\rm ht}_{\overline{\mathcal L}}(T(m)) + c \int_X \widetilde{\Xi }(v,z,m))d\mu \\\nonumber
&= c \sigma _1(m)\pi ^2(2/3)(\Ei(-4\pi |m|v)+e^{-4\pi |m|v}/(4\pi |m|v)).
\end{align*}
In both cases we get the equality we claimed with $c = 18/\pi ^2.$\hfill$\Box$\\

\nn {\bf Remark.}\label{Remvolhyp} The constant $c$ is explained by the fact that in the context of \cite{BKK} and \cite{BBK} one has
\begin{align*} 
\int c_1( \overline{\mathcal L})^2 = c_1( \overline{\mathcal L})\cdot c_1( \overline{\mathcal L}) = T(1)\cdot T(1) = 2
\end{align*} 
with
\begin{align*} 
c_1( \overline{\mathcal L}) = 12 \,(dx_1dy_1/(4\pi y_1^2) + dx_2dy_2/(4\pi y_2^2))
\end{align*} 
and
\begin{align*} 
{c}_1(\overline{\mathcal L})^2 = 2\cdot (12/4\pi )^2 (dx_1dy_1 dx_2dy_2/(y_1 y_2)^2) = (18/\pi ^2 )d\mu (z).
\end{align*} 

Finally, we have all the material for the\\

{\bf Proof of the Main Theorem.} As we know the modularity of both sides of \eqref {eq:kudla-modi} the equality 
for $m\not=0$ following from Theorem \ref {CorKucon1} is sufficient (and, hence, gives the value for $m=0$ we up to now did not 
determine directly).\hfill$\Box$\\

\section{Remarks towards a direct calculation of the constant term}

Though we don't really need this, to strive for some completeness, we will make some remarks concerning the case $m=0.$\\

For $m=0,$ as a consequence of the log-log-singularity  of the metric on $\overline{\mathcal L},\, {\rm ht}_{\overline{\mathcal L}}(T(0))$ 
is not defined and by the same reason $\int_X \widetilde{\Xi}_\rho (v,z,0)d\mu$ does not exist.
 Therefore instead of \eqref{cxcxZm},  we have to use the formula   
 \begin{align*}\label{cxcxZ0}
\widehat{c}_1( \overline{\mathcal L})^2\cdot  \widehat{Z}_\rho (0) 
=&\int_X (\widetilde{\Xi}_\rho (v,z,0)+  (1/24 - 1/(8\pi v)) g_0)d\mu \\\nonumber
&+ \widetilde{c} (\zeta '(-1)/\zeta (-1)+1/2),\nne
\end{align*}
where
\begin{align*}\label{eqg0}
g_0 :=  -\log \parallel \Delta (z_1)\Delta (z_2)\parallel ^2.\nne
\end{align*} 
This formula comes out as here one has
\begin{align*}
\widehat{Z}_\rho (0) &= (T(0), \widetilde{\Xi}_\rho )\\
\widehat{T}(0) &= (T(0), c_0g_0),\,\,c_0 := -1/24+1/8\pi v
\end{align*} 
and
 \begin{align*}\label{cxZ00}
\widehat{c}_1( \overline{\mathcal L})^2\cdot  \widehat{Z}_\rho (0) 
&= \widehat{c}_1( \overline{\mathcal L})^2\cdot (\widehat{T}(0) + (\widehat{Z}_\rho (0) - \widehat{T}(0))\\\nonumber
&= \widehat{c}_1( \overline{\mathcal L})^2\cdot \widehat{T}(0) +  \widehat{c}_1( \overline{\mathcal L})^2\cdot (\widehat{Z}_\rho (0) - (\widehat{T}(0))\nne
\end{align*}
where the first summand is known to be (\cite{BKK} Theorem 7.61) 
 \begin{align*}
 \widehat{c}_1( \overline{\mathcal L})^2\cdot \widehat{T}(0) &= \widehat{c}_1( \overline{\mathcal L})^2\cdot c_0\widehat{c}_1( \overline{\mathcal L})\\
 &= \widetilde{c}(\zeta '(-1)/\zeta (-1)+1/2),\,\widetilde{c}= -(1/2)12^2c_0
\end{align*}
and (using the formulae (3.3.1), (3.1.1), and (3.0.5) from Part I) the last one gives the integral above in the formula (\ref{cxcxZ0})
 \begin{align*}\label{cxcxZ00}
 \widehat{c}_1( \overline{\mathcal L})^2\cdot (\widehat{Z}_\rho (0) - (\widehat{T}(0))
&= \int_X(\widetilde{\Xi}_\rho (v,z,0)-c_0g_0)c_1( \overline{\mathcal L})^2\\\nonumber
&= \int_X(\Xi (v,z,0) -c_0g_0)c_1( \overline{\mathcal L})^2 + \int_X \rho (z)\check{\Xi }^+ (v,z,0)c_1( \overline{\mathcal L})^2\nne
\end{align*}

We already fixed 
\begin{align*}
f_\rho(0) = \int_X \rho (z)\check{\Xi }^+ (v,z,0))d\mu.
\end{align*}
Thus, if one wants to avoid the reasoning from the proof above, for a direct proof of the $m=0-$case 
it remains to show that $(\pi ^2/18)A'(v,z,0)$ has the same value as
\begin{align*} 
\label{cxcxZ}
\int_X ({\Xi} (v,z,0) &+ (1/24 - 1/(8\pi v)) g_0)d\mu  \\\nonumber
&= \int_X (({\Xi} ^*(v,z,0)-(1/24)g_0) + (1/24 - 1/(8\pi v)) g_0)d\mu  \\\nonumber
&= \int_X ({\Xi} ^*(v,z,0) - 1/(8\pi v) g_0)d\mu  . \nne
\end{align*} 
Observe, if we split the integrals, then we would get two divergent integrals
where for the integral over $g_0$ the relevant terms in an asymptotic expansion at 
the boundary had been calculated already in \cite{K}. \\

\begin{appendix}
 
\setcounter{section}{0}

\section{Appendix: Fourier expansion of $E_2(\tau ,s)$}
Now, here we add the direct proof of Theorem \ref{ThderEiss} where in the main part we relied on Maple calculations.\\

\nn {\bf Theorem.}  \emph{We have
\begin{align*}
E_2(\tau ,1) = \sum_{m\in\Z}a(v,1,m) q^m\nne
\end{align*}
and, denoting by $E'_2(\tau ,s)$ the derivative of $E_2(\tau ,s)$ with respect to $s,$ we get
\begin{align*}
E'_2(\tau ,1) = \sum_{m\in\Z}a'(v,1,m) q^m\nne
\end{align*}
with
\begin{align*}
a(v,1,m) &=
\begin{cases} \sigma _1(m)\,\,&{\rm for} \,\,m>0\\\nonumber
- 1/24 + 1/(8\pi v)&{\rm for} \,\,m=0\\\nonumber
 0\,\,&{\rm for}\,\,m<0\nonumber
 \end{cases}\\
a'(v,1,m) &= \begin{cases} \sigma _1(m)(1/(4\pi mv) + \sigma '/\sigma ) \,\,&{\rm for} \,\,m>0\\\nonumber
 -(1/24)(24\zeta '(-1)+\gamma -1+\log(4\pi v))\\\nonumber
 \quad \,-(1/(8\pi v))(-\gamma +\log(4\pi v))&{\rm for} \,\,m=0\\\nonumber
 \sigma _1(m)(\Ei(-4\pi |m|v) + 1/(4\pi |m|v)e^{-4\pi |m|v})\,\,&{\rm for}\,\,m<0\nne
 \end{cases}
\end{align*}
where with $\sigma ^\ast $ as in \eqref {sigast}}
\begin{align*}\label{eqdersigma}
\sigma := \sigma ^\ast _{1/2}(m),\,\,\sigma ' := {\sigma ^\ast }'_{1/2}(m).\nne
\end{align*}\\


{\bf Proof. 1.}
We have (see for instance Iwaniec \cite{Iw} p.205)
\begin{align*}
K_\nu  (t) :=& \int_0^\infty e^{-t\cosh u} \cosh (\nu u) du \\[.3cm]
=& \frac{\sqrt \pi (t/2)^\nu }{\Gamma (\nu +(1/2))} \int_1^\infty e^{-tr} (r^2-1)^{ \nu -(1/2)}dr .
\end{align*}
Hence, from(\ref{EisFourdev}) we get
\begin{align*}
E^\ast (\tau ,s) = & v^s \zeta ^\ast (2s) + v^{1-s} \zeta ^\ast (2s-1) \\[.3cm]
& + \sum _{m\in\Z, m\not=0} 2\sigma ^\ast _{s-(1/2)} (|m|) \frac{(v\mid m \mid \pi)^s}{\Gamma (s)\sqrt{|m| }}
 \int_1^\infty e^{-2\pi\mid m \mid vr} (r^2-1)^{s - 1}dr \,e^{2\pi imu}.
 \end{align*}
As in the rudimentary proof of Theorem \ref{ThderEiss} we abbreviate
\begin{align*}
c_0(v,s) :=& \,\,v^s \zeta ^\ast (2s) + v^{1-s} \zeta ^\ast (2s-1) ,\\
\check c_0(v,s) :=& \,\,\partial _vc_0(v,s) = sv^{s-1} \zeta ^\ast (2s) + (1-s)v^{-s} \zeta ^\ast (2s-1) ,\\
 \end{align*}
 and, for $m\not=0,$
\begin{align*}
\alpha :=& \,2\pi |m|v,\\
c_m(v,s) :=& \,2\sigma ^\ast _{s-(1/2)} (|m|) \frac{(v|m|\pi)^s}{\Gamma (s)\sqrt{|m| }},\\
I_m (v,s):=&  \int_1^\infty e^{-2\pi\mid m \mid vr} (r^2-1)^{s - 1} dr ,\\
J_m (v,s):=&  \int_1^\infty e^{-2\pi\mid m \mid vr} (r^2-1)^{s - 1}r dr ,
 \end{align*}
 
 and get
 \begin{align*}
 E_2(\tau ,s) = - (1/(4\pi))( \check c_0(v,s) + {\sum }'_m ((s/v \,\,-& 2\pi m)c_m(v,s)I_m(v,s) \\[.3cm]
-& 2\pi |m| c_m(v,s)J_m(v,s))e(mu)
 \end{align*}
  i.e.
 \begin{align*}\label{E2s}
 E_2(\tau ,s) =& - \frac{1}{4\pi }(\check c_0(v,s)\\\nonumber
&+ \sum_{m>0} ((s/v)c_m(v,s)I_m(v,s) - 2\pi |m| c_m(v,s)(I_m(v,s)+J_m(v,s))e(mu)\\\nonumber
&+ \sum_{m<0} ((s/v)c_m(v,s)I_m(v,s) + 2\pi |m| c_m(v,s)(I_m(v,s)-J_m(v,s))e(mu)\,).\nne
 \end{align*}
 
{\bf 2.} Denoting by ' the derivation with respect to $s,$ one has
 \begin{align*}\label{E'2s}
 E'_2(\tau ,s) = - \frac{1}{4\pi }(\!\! &\,\,\check c'_0(v,s)&&\\[.3cm]
&+\,\,\sum_{m>0} (((1/v)c_m &+& (s/v)c'_m)I_m - 2\pi |m| c'_m(I_m+J_m)\\
&&+& (s/v)c_mI'_m - 2\pi |m| c_m(I'_m+J'_m))e(mu)\\[.3cm]
&+\,\, \sum_{m<0} (((1/v)c_m &+& (s/v)c'_m)I_m + 2\pi |m| c'_m(I_m-J_m)\\
&&+& (s/v)c_mI'_m + 2\pi |m| c_m(I'_m-J'_m))e(mu)\,).\nne 
 \end{align*}
 
From
$$
c_m(v,s) := 2\sigma ^\ast _{s-(1/2)} (\mid m \mid ) \frac{(v\mid m \mid \pi)^s}{\Gamma (s)\sqrt{|m| }}
$$
we get
$$
c_m(v,1) = \,\sigma ^\ast _{1/2} (\mid m \mid ) \alpha /\sqrt{|m| }
$$
and with $\sigma _s := \sigma ^\ast _s(|m| )$
$$
c'_m(v,s) = ((\sigma '_{s-(1/2)}/\sigma _{s-1/2}) + \log (\alpha /2) - (\Gamma '(s)/\Gamma (s)))c_m(v,s).
$$
Using $\Gamma '(1) = - \gamma , \gamma $ the Euler constant, we have
$$
c'_m(v,1) = ((\sigma '_{1/2}/\sigma _{1/2}) + \log (\alpha /2) + \gamma )c_m(v,1).
$$
Now, we can write for $s = 1$
$$
\begin{array}{rrl}
 E'_2(\tau ,1) = &- (1/(4\pi))( \check c'_0(v,1)&\\[.3cm]
&+\sum_{m>0} ((1/v)c_m ((1 +& (\sigma '_{1/2}/\sigma _{1/2}) + \log (\alpha /2) + \gamma )I_m\\
& -& \alpha ((\sigma '_{1/2}/\sigma _{1/2}) + \log (\alpha /2) + \gamma )(I_m+J_m)\\
&+& I'_m - \alpha (I'_m+J'_m))e(mu)\\[.3cm]
&+ \sum_{m<0} ((1/v)c_m ((1 +& (\sigma '_{1/2}/\sigma _{1/2}) + \log (\alpha /2) + \gamma )I_m\\
& +& \alpha ((\sigma '_{1/2}/\sigma _{1/2}) + \log (\alpha /2) + \gamma )(I_m-J_m)\\
&+ &I'_m + \alpha (I'_m-J'_m))e(mu)\,).\\[.3cm] 
\end{array}
$$
One has to determine the values at $s=1$ of the functions in this relation: 
From the definitions one has 
$$
\begin{array}{rl}
I_m(v,1) =& e^{-\alpha } (1/\alpha ),\\
J_m(v,1) =& e^{-\alpha } ((1/\alpha ) + (1/\alpha ^2)),
\end{array}
$$
i.e.
$$
\begin{array}{rl}
I_m(v,1) - J_m(v,1) =& - e^{-\alpha } (1/\alpha ^2),\\
I_m(v,1) + J_m(v,1)  =& e^{-\alpha } ((2/\alpha ) + (1/\alpha ^2)),
\end{array}
$$
and, hence, via \eqref{E2s} immediately the formulae for the $a(v,1,m), m\not=0,$ in the theorem. \\

{\bf 3.} For the other terms we will use the well known relation
$$
\Gamma '(1) = \int_0^\infty e^{-t} \log t \,dt = -\gamma 
$$
and its consequence
$$
\int_0^\infty e^{-\alpha t} \log t \,dt = - (1/\alpha )(\gamma + \log \alpha ).
$$
Moreover, one has using partial integration
$$
\begin{array}{rl}
\int_0^\infty e^{-\alpha t} t\log t \,dt =&  (1/\alpha )((1/\alpha ) - (1/\alpha )(\gamma + \log \alpha )),\\
\end{array}
$$
and with $ - {\rm Ei} (-s) = \int_1^\infty e^{-st} dt/t$
$$
\begin{array}{rl}
\int_2^\infty e^{-\alpha t} \log t \,dt =& (1/\alpha )( - {\rm Ei} (-2\alpha ) + e^{-2\alpha }\,\log 2 )\\[.4cm]
\int_2^\infty e^{-\alpha t} t\log t \,dt 
=& (1/\alpha ^2)( - {\rm Ei} (-2\alpha ) + e^{-2\alpha }(2\alpha \log 2 + 1 + \log 2) ).\\
\end{array}
$$
Hence we get
$$
\begin{array}{rl}
I'_m(v,s) =& \int_1^\infty e^{-\alpha r} (r^2-1)^{s-1} \log(r^2-1) \,dr\\[.4cm]
I'_m(v,1) =& \int_1^\infty e^{-\alpha r} \log(r-1) \,dr + \int_1^\infty e^{-\alpha r} \log(r+1) \,dr\\[.4cm]
=& e^{-\alpha }\int_0^\infty e^{-\alpha r} \log r \,dr +e^{\alpha }\int_2^\infty e^{-\alpha r} \log r \,dr\\[.4cm]
=& e^\alpha (1/\alpha )( - {\rm Ei} (-2\alpha )) - e^{-\alpha }(1/\alpha )(\log (\alpha /2) + \gamma )
\end{array}
$$
and similarly
$$
\begin{array}{rl}
J'_m(v,s) =& \int_1^\infty e^{-\alpha r} (r^2-1)^{s-1} r\log(r^2-1) \,dr\\[.4cm]
J'_m(v,1) =& \int_1^\infty e^{-\alpha r} r\log(r-1) \,dr + \int_1^\infty e^{-\alpha r} r\log(r+1) \,dr\\[.4cm]
=& e^{-\alpha }\int_0^\infty e^{-\alpha r} (r+1)\log r \,dr +e^{\alpha }\int_2^\infty e^{-\alpha r} (r-1)\log r \,dr\\[.4cm]
=& e^\alpha ({\rm Ei} (-2\alpha ))((1/\alpha )-(1/\alpha ^2) )
+ e^{-\alpha }(1/\alpha ^2)(2 - (\log (\alpha /2) + \gamma )(1+\alpha )).
\end{array}
$$
One has
$$
\begin{array}{rl}
I'_m(v,1)+J'_m(v,1) =& e^\alpha (-{\rm Ei}(-2\alpha ) (1/\alpha ^2) )\\
&+ e^{-\alpha }(1/\alpha ^2)(2 - (\log (\alpha /2) + \gamma ) - 2\alpha (\log (\alpha /2) + \gamma ))\\[.3cm]
I'_m(v,1)-J'_m(v,1) =& e^\alpha ({\rm Ei} (-(2\alpha ))((1/\alpha ^2)-(2/\alpha ) )\\
&- e^{-\alpha }(1/\alpha ^2)(2 - (\log (\alpha /2) + \gamma )).
\end{array}
$$
Hence, from \eqref{E'2s}, in the equation for $E'_2(\tau ,1)$ the coefficient of $c_m(v,1)/v$ for $m>0$ comes out as
$$
e^{-\alpha }(-2\sigma '/\sigma -(1/\alpha ))
$$
and for $m<0$ as
$$
e^{-\alpha }(-(1/\alpha )) - e^\alpha 2{\rm Ei} (-(2\alpha )).
$$
Remembering that for $m>0$ with $q = e(u+iv)$ one has $q^m = e^{-\alpha }e(mu)$ and 
for $m<0$ $q^m = e^\alpha e(mu),$ we get with $\sigma := \sigma ^\ast _{1/2}(|m| )$
$$
\begin{array}{rl}
E'_2(\tau ,1) =& -(1/4\pi )\check c'_0(v,1)\\[.3cm]
&+ \sum _{m>0}(\sqrt{|m|}/(4\pi ))(\sigma /(4\pi mv) + \sigma ') q^m\\[.3cm]
&+ \sum _{m<0}(\sqrt{|m|}/(2\pi ))\sigma (-{\rm Ei} (-(2\alpha )) + e^{4\pi mv} /(4\pi mv) ) q^m.\\
\end{array}
$$
One can simplify this a bit using the relation $\sigma = \sigma ^\ast _{1/2} = \sigma ^\ast _{-1/2}$ and, with
$$
\sqrt {|m|} \sigma = \sigma _1 = \sum _{d\mid m} d.
$$
have our claim for $m\not=0.$\\
{\bf 4.} We still have to treat the case $m=0,$ i.e., starting with 
\begin{align*}
c_0(v,s) :=& v^s \zeta ^\ast (2s) + v^{1-s} \zeta ^\ast (2s-1) ,\\[.3cm]
\check c_0(v,s) :=& \partial _vc_0(v,s) = sv^{s-1} \zeta ^\ast (2s) + (1-s)v^{-s} \zeta ^\ast (2s-1) ,
 \end{align*}
determine $\check c_0(v,1)$ and $\check c'_0(v,1).$\\

 We look at
\begin{align*}
a_1(s) := sv^{s-1}\zeta ^\ast (2s) = \zeta (1-2s)\Gamma (1/2-s)\pi^{s-1/2}sv^{s-1}
\end{align*}
and get, using standard material assembled in the Zeta Tool Remarks below,
\begin{align*}
a_1(1) &= - \zeta (-1)2\pi = \pi/6,\\
a'_1(s) &= (-2\zeta '(1-2s)/\zeta (1-2s) - \Gamma '(1/2-s)/\Gamma (1/2-s) + \log \pi + 1/s + \log v)a_1(s)\\
a'_1(1) &= (-2\zeta '(-1)/\zeta (-1) - \Gamma '(-1/2)/\Gamma (-1/2) + 1 + \log (\pi v))a_1(1)\\
 &= (\pi /6)(24\zeta '(-1) + \gamma  - 1 + \log (4\pi v)).
\end{align*}
Similarly, we take
\begin{align*}
a_2(s) := (1-s)v^{-s}\zeta ^\ast (2s-1) = (1-s)\zeta (2s-1)\Gamma (s-1/2)\pi^{1/2-s}v^{-s}
\end{align*}
and for
\begin{align*}
F(s) := (1-s)\zeta (2s-1) = (1-s)(1/(2(s-1)) + \gamma +\dots)
\end{align*}
get
\begin{align*}
F(1) := - 1/2,\,\,F'(1) = - \gamma ,
\end{align*}
while for
\begin{align*}
G(s) := \Gamma (s-1/2)\pi^{1/2-s}v^{-s}
\end{align*}
one has
\begin{align*}
G(1) &= \Gamma (1/2)\pi ^{-1/2}v^{-1} = 1/v\\
G'(s) &= (\Gamma '(s-1/2)/\Gamma (s-1/2)-\log \pi -\log v)G(s)\\
G'(1) &=  (\Gamma '(1/2)/\Gamma (1/2)-\log (\pi v))1/v\\
 &=  (1/v)(-\gamma - \log (4\pi v)),
\end{align*}
and, hence,
\begin{align*}
a_2(1) &= F(1)G(1) = -1/(2v) = \zeta (-1)4\pi (3/(2\pi v))\\
a_2'(1) &= F'(1)G(1) + F(1)G'(1) = -\gamma (1/v) -(1/2)((-\gamma -  \log (4\pi v))1/v\\
&= (1/2v)(-\gamma + \log (4\pi v)).
\end{align*} 
Finally we have the claim for the constant terms
\begin{align*}
\check{c}_0(v,1) &= a_1(1)+a_2(1) = \zeta (-1)(-2\pi + 6/v) = \pi /6 - (1/2v)\\ 
\check{c}'_0(v,1) &= a'_1(1)+a'_2(1) \\
&= (\pi /6)(24\zeta '(-1) + \gamma  - 1 + \log (4\pi v)) + (1/2v)(-\gamma + \log (4\pi v)).\\
\end{align*}

\hfill$\Box$\\
\end{appendix}

\end{document}